\documentclass[12pt]{amsart}

\usepackage{color}
\definecolor{darkgreen}{rgb}{0, 0.40, 0}

\usepackage[
ps2pdf,                
colorlinks=true,
linkcolor=blue,
citecolor=darkgreen,
urlcolor=magenta
]{hyperref}

\usepackage{epsfig}
\usepackage{latexsym}
\usepackage{psfrag}
\usepackage{amssymb}
\usepackage{amsmath}
\usepackage{amsthm} 
\usepackage{verbatim} 

\newcommand{\calA}{{\mathcal{A}}}

\newcommand{\calC}{{\mathcal{C}}}
\newcommand{\calD}{{\mathcal{D}}}

\newcommand{\calM}{{\mathcal{M}}}

\newcommand{\calP}{{\mathcal{P}}}
\newcommand{\calQ}{{\mathcal{Q}}}

\renewcommand{\setminus}{{\smallsetminus}}

\newcommand{\bigmid}{{\:\mid\:}}
\newcommand{\st}{{\bigmid}}
\newcommand{\from}{{\colon}} 

\newcommand{\nin}{{\notin}}

\newcommand{\bdy}{{\partial}} 

\newcommand{\diam}{{\operatorname{diam}}} 

\newcommand{\PML}{{\mathcal{PML}}}

\newsavebox{\savepar}

\theoremstyle{plain}
\newtheorem{theorem}{Theorem}[section]

\newtheorem{lemma}[theorem]{Lemma}

\theoremstyle{definition}

\newtheorem{remark}[theorem]{Remark}

\newcommand{\fakethm}[2]
{
\vspace{6 pt}
\noindent
{{\bf Theorem~\ref{#1}.}~{\em #2}}
\vspace{6 pt}
}

\usepackage{psfig}
\def\givenname#1{}
\def\surname#1{}

\long\def\realfig#1#2#3{
\begin{figure}[htbp]
\centerline{\psfig{file=Figures/#1.ps}}
\caption[#1]{#3}
\label{#2}
\end{figure}}

\newcommand{\ovX}{{\overline{X}}}
\newcommand{\ovY}{{\overline{Y}}}
\newcommand{\boundary}{\partial}
\newcommand{\intersect}{\cap}
\newcommand{\union}{\cup}

\begin{document}

\title{High distance knots}

 \author{Yair N. Minsky}
\givenname{Yair}
\surname{Minsky}
 \address{Department of Mathematics\\
       Yale University\\
       New Haven, CT 06520-8283}
 \email{yair.minsky@yale.edu}

 \author{Yoav Moriah}
\givenname{Yoav}
\surname{Moriah}
 \address{Department of Mathematics\\
       Technion\\
       Haifa 32000 Israel}
 \email{ymoriah@tx.technion.ac.il}

 \author{Saul Schleimer}
\givenname{Saul}
\surname{Schleimer}
 \address{Mathematics Institute\\
Zeeman Building\\
University of Warwick\\
Coventry CV4 7AL UK}
 \email{S.Schleimer@warwick.ac.uk}

\date{\today}

\begin{abstract} We construct knots in $S^3$ with Heegaard splittings
  of arbitrarily high  distance, in any genus. As an application, for
  any positive integers $t$ and $b$ we find a  tunnel number $t$ knot
  in the three-sphere which has no $(t,b)$--decomposition.
\end{abstract}


\maketitle

\section{Introduction}
\label{Sec:Intro}

In this paper we address the problem of generating knots in $S^3$ with
high distance, in the sense of Hempel. 
Without any restriction on the ambient manifold, knots with high
distance are easy to construct: 
Hempel~\cite{Hempel01}, adapting an idea of
Kobayashi~\cite{Kobayashi88b}, constructs Heegaard splittings $(V,W)$
of arbitrarily high distance. Given such a splitting remove an
unknotted solid torus from, say, $V$ to obtain a compression body
$V_0$.  The result is a knot space in some manifold $M$ and clearly
the distance of the Heegaard splitting $(V_0,W)$ is at least the
distance of $(V,W)$.  See Section~\ref{Sec:Preliminaries} for
definitions.

The problem becomes more challenging if the ambient manifold $M$ is
specified beforehand.  The argument above fails if
$M$ does not admit any high distance splitting.  The case of
$S^3$ is of particular interest, and here we know \cite{Waldhausen68}
that every Heegaard
splitting is isotopic to a standard one.  The disk complexes
for a standard splitting have distance zero, and indeed an infinite-diameter
intersection. This makes finding an appropriate compression body more
difficult. Nevertheless, we prove:

\fakethm{Thm:HighDistanceKnots}{For any pair of integers $g > 1$ and
  $n > 0$ there is a knot $K \subset S^3$ and a genus $g$ splitting of
  its exterior $E(K)$ having distance greater than $n$.}

We can consider another measure of complexity of a knot: $K \subset M$
has a {\em $(g, b)$--decomposition} if $K$ may be isotoped to have
exactly $b$ bridges with respect to a genus $g$ Heegaard splitting of
$M$.  When $b = 0$ we further insist that $K$ be a {\em core} of one
of the handlebodies (see Remark~\ref{Rmk:g0}). The classical notion of
``bridge position'' corresponds to a $(0,b)$--decomposition.

Work of Scharlemann-Tomova \cite{ScharlemannTomova05} and Tomova 
\cite{Tomova05,Tomova07} links these two notions of complexity. It implies
in particular that the existence of a $(g,b)$--decomposition gives
upper bounds on the distance of Heegaard splittings of genus greater
than $g$ (see Section~\ref{Sec:Decompositions} for details).  With
this we obtain the following consequence of
Theorem~\ref{Thm:HighDistanceKnots}:

\fakethm{Thm:GBKnots}{For any positive integers $t$ and $b$ there is a
knot $K \subset S^3$ with tunnel number $t$ so that $K$ has no
$(t, b)$--decomposition.}

This theorem answers question 1.9 of Kobayashi-Rieck \cite{KobayashiReick04b}.
Note that any knot $K$ with tunnel number $t $ (see
Section~\ref{Sec:Preliminaries}) has by definition a $(t+1,
0)$--decomposition. The corresponding splitting surface of $E(K)$
has minimal genus.

Previous constructions have yielded tunnel number $t$ knots without
$(t,1)$--decompositions (Moriah and
Rubinstein~\cite{MoriahRubinstein97}, Morimoto, Sakuma, and
Yokota~\cite{MSY96}, and Eudave-Mu\~noz~\cite{Eudave-Munoz02}).  The
knots described in~\cite{MoriahRubinstein97} have
$(t,2)$--decompositions.  The knots in~\cite{MSY96} are all tunnel
number one knots and have $(1,2)$--decompositions.  The knots
in~\cite{Eudave-Munoz02} are all tunnel number one and some are known
to be $(1,2)$ while the rest have unknown optimal
$(1,b)$--decomposition.

More recently Eudave-Mu\~noz \cite{Eudave-Munoz07} has exhibited tunnel number
one knots in $S^3$ which are not $(1,2)$.  His examples are either $(1,3)$ or
$(1,4)$ knots but exactly which is not yet known.

Johnson and Thompson have shown, in~\cite{JohnsonThompson06} (see also
Johnson~\cite{Johnson06}), that for arbitrarily large $n$ there are
tunnel number one knots which have no $(1,n)$--decomposition.  Their
proof also uses the Scharlemann-Tomova results to relate Heegaard
distance to bridge position. Constructing high distance 
tunnel number one knots is 
simplified by the 
special properties of 
disk complexes in genus 2 compression bodies -- see the discussion in
\S\ref{Sec:HighDistance}.

\subsection*{Further questions}

One is tempted to speculate that Theorem~\ref{Thm:HighDistanceKnots}
holds for any 3-manifold. We note however that our construction
fails for the standard Heegaard splitting of a connected sum of
$S^2\times S^1$'s, because the disk complexes on
both sides of that splitting are identical. A more plausible
conjecture is that, after stabilizing once or possibly twice, any
Heegaard splitting will admit a construction like the one we have
used. 

\subsection*{Acknowledgements}
This work arose out of conversations between the authors at the
Workshop on Heegaard Splittings of $3$-Manifolds, held at the Technion
in Haifa (Israel), July of 2005.  Further work was done at Rutgers and
Yale University during February of 2006.  We wish to thank the
Technion, Rutgers, and Yale University for their hospitality.
This work was partially supported by NSF grants DMS-0504019
and DMS-0508971, and Grant No.~2002039 from the
US-Israel Binational Science Foundation (BSF).
It is in the public domain.

\section{Preliminaries} 
\label{Sec:Preliminaries}

A {\em compression body} is a $3$-manifold $V$ obtained from a
surface $S$ cross an interval $[0,1]$ by attaching a finite number of
2-handles and 3-handles to $S\times\{0\}$.
The component $S\times\{1\}$ of the boundary will be denoted by
$\partial_{+}V$, and 
$\partial V \setminus \partial_{+}V$ will be denoted by  $\partial_{-}V$. 
The trivial cases
where $V$ is a handlebody or $V = S \times [0,1]$ are allowed.

A {\em Heegaard splitting} for a $3$-manifold $M$ is a decomposition
$M = V \cup_S W$ where $V, W$ are compression bodies so that $S =
\bdy_+ V = \bdy_+ W = V \cap W$.  The surface $S$ will be called the
{\em Heegaard surface} of the Heegaard splitting.

\subsection{Complexes}

Let $C_S$ be the 1-skeleton of Harvey's {\em complex of curves}
(see~\cite{Harvey81}). That is, given a closed connected orientable
surface $S$ of genus at least two, let $\calC_S$ be the graph whose
vertices are isotopy classes of essential simple closed curves and
whose edges connect distinct vertices with disjoint representatives.
As usual, we only distinguish curves from their isotopy classes when
necessary. 
 
We remark that $\calC_S$ is connected.  Place a metric $d( \cdot,
\cdot)$ on $\calC_S$ by setting the length of every edge to be one.
For subsets $X, Y \subset \calC_S$ we define $d(X, Y) = \min \{ d(x,
y) \st x \in X,~y \in Y \}$.

An essential curve $\alpha \subset S = \bdy_+ V$ is a {\it meridian}
of the compression body $V$ if $\alpha$ bounds an essential disk in
$V$.  Given a compression body $V$ so that $\partial_+V = S$, let
$\calD_V$ be the sub-complex of $\calC_S$ spanned by meridians of $V$.
This is the one-skeleton of McCullough's {\em disk complex}
(see~\cite{McCullough91}).

Hempel~\cite{Hempel01} defined the {\em distance} of a
Heegaard splitting $V \cup_S W$ to be
$$d(V, W) = d(\calD_V, \calD_W).$$

\subsection{Laminations}

We let $\PML(S)$ denote the space of projective measured laminations
on S. For a good reference on measured laminations see
Bonahon~\cite{Bonahon01}, Casson and Bleiler~\cite{CassonBleiler88},
or Hatcher~\cite{Hatcher88}.

We will need the following lemma due to Hempel~\cite{Hempel01}, who
generalized an argument of Kobayashi~\cite{Kobayashi88b}.

\begin{lemma}
\label{Lem:Hempel}
Suppose that $X, Y \subset \calC_S$,  and let $\ovX$ and $ \ovY$
denote their 
closures in $\PML(S)$.  Let  $\Phi$ be a pseudo-Anosov map with stable and
unstable laminations $\lambda^\pm$.  Assume that $\lambda^- \nin~\ovY$
and $\lambda^+ \nin~\ovX$.  Then $d(X, \Phi^n(Y)) \to \infty$ as $n
\to \infty$.  
\qed
\end{lemma}

This statement is virtually identical to that of Theorem~2.4 in
Abrams-Schleimer~\cite{AbramsSchleimer02}, the only difference being
that $X=Y$ there. The same proof carries through.

\subsection{Pants, waves and seams} 
Given a surface $S$ and an essential subsurface $Y\subset S$, a {\em
  wave} in $Y$ is a properly embedded arc in $Y$ not homotopic rel
  $\boundary Y$   into $\boundary Y$, and 
with both ends incident to the same boundary component of $Y$ from the
  same side.

A {\em pair of pants} is a three-holed sphere. 
A {\em pants decomposition} $\calP = \{ \alpha_i \}$ of a surface $S$
is a maximal collection 
of  essential simple closed curves which are disjoint and not
parallel; hence (we are assuming genus at least two) $\overline{S -
  \calP}$ is a union of pairs of pants.  

A {\em seam} in a pair
of pants is an essential properly embedded arc connecting
distinct boundary components.
Up to isotopy rel
boundary, a pair of pants has three distinct waves and three distinct
seams. 

We say an essential curve or a lamination $\alpha\subset S$ {\em
  traverses} or {\em 
  has} a wave (seam) with respect to $Y$  if $\alpha$ intersects $Y$
  minimally in its isotopy class and a component of $\alpha\intersect
  Y$ is a wave (seam).

A standard outermost bigon argument proves the existence of waves
among meridians in a compression body:

\begin{lemma}
\label{Lem:Wave}
Let $V$ be a compression body and
let
$\calM$ be a maximal collection of non-parallel, disjoint
meridians on $S=\partial_+ V$. 
If $\alpha$ is a meridian of $V$ then $\alpha$ is either parallel to a
component
of $\calM$  
or $\alpha$ has a wave with respect to one of the pants components of
$S\setminus \calM$. 
\qed
\end{lemma}

\noindent Let $V, S$ and $\calM$ be as above, and let 
$\overline{\calD_V}$ denote closure in $\PML(S)$. 
Here and in what follows, by a ``seam of a curve
system'' we mean a seam of a pants component of its complement.
We then have: 

\begin{lemma}
\label{Cor:FullType}
No lamination $\lambda \in \overline{\calD_V}$ traverses all seams of
$\calM$.
\end{lemma}

\begin{proof}
Suppose that $\{\beta_i\} \in \calD_V$ is a sequence of meridians
converging to $\lambda \in \PML(S)$.  If $\lambda \in \calM$ then
clearly $\lambda$ does not traverse all seams of
$\calM$. 
By passing to a subsequence, we may assume that none of the $\beta_i$
lie in $\calM$.  It follows from Lemma~\ref{Lem:Wave} that all
$\beta_i$ have a wave in some pants of $S \setminus \calM$.  Passing
to subsequences again we may assume that all of the $\beta_i$ have the
same wave $w$ in the same pair of pants $Y$. 

If $\lambda$ traversed all seams
of $\calM$, then
for large enough index $i$ the curve
$\beta_{i}$ would also 
traverse all these seams. Since one of
the seams of $Y$ intersects $w$ essentially it would follow that
$\beta_{i}$ self-intersects, a contradiction.
\end{proof}

\subsection{$(g,b)$--Decompositions}

Suppose $\calA \subset V$ is a disjoint collection of properly
embedded arcs in a compression body $V$ where $\bdy \calA \subset
\bdy_+ V$.  We say $\calA$ is {\em unknotted } if $\calA$ can be
properly isotoped, rel boundary, into $\bdy V$.

The following well-known generalization of bridge position is
due to Doll~\cite{Doll92}: 
Suppose that $M = V \cup_S W$ and $K$ is a knot in $M$.  The knot $K$
is in {\em bridge position with respect to $S$} if $K$ is transverse
to $S$ and either
\begin{itemize}
\item[{\em (i)}] $K \cap S \neq \emptyset$, and both $K \cap V$ and $K
\cap W$ are unknotted, or
\item[{\em (ii)}] $K \cap S = \emptyset$, 
$K \subset V$ (without loss of generality), and $V\setminus n(K)$ is a
  compression  body. 
\end{itemize}
If $g = g(S)$ and $b = |K \cap S|/2$ then we say that $K$ admits a
  {\em $(g, b)$--decomposition}.

\begin{remark}\rm 
Our definition of a $(g,0)$--decomposition
  is non-standard, due to the additional requirement on $V\setminus n(K)$. In
  particular this definition makes a $(g,0)$--decomposition equivalent
  to a Heegaard splitting of $E(K)$ of genus $g$. 
\label{Rmk:g0}
\end{remark}

\subsection{Tunnel number}

Given a knot $K$ in a $3$-manifold $M$, a {\it tunnel system} for $K$
is a collection $\mathcal{T } = \{a_1, . . . , a_n\}$ of disjoint
arcs, properly embedded in $M\setminus K$, 
such that $M \setminus n(\mathcal{T} \cup K)$ is a handlebody. 
The {\em tunnel number} $t(K)$ of the knot is the minimal cardinality
$n$ of such a tunnel system.

We note that a knot with tunnel number $t$ has a $(t+1,0)$--decomposition, and
that $t+1$ is the minimal genus for which this is 
possible. 
Hence it becomes interesting to reduce the genus by one, and ask what
$(t,b)$--decompositions such a knot admits. 

Of course, if $K$ has a $(g, b)$--decomposition then $K$ also has $(g,
b+1)$ and $(g+1, b)$--decompositions.  With a bit more care 
one can check that any knot with a $(g, b)$--decomposition also has
a $(g+1, b-1)$--decomposition. 
On the other hand, when going from a $(g,b)$ to a $(g-1,b')$--decomposition,
$b'$ may need to grow arbitarily. This follows from 
Theorem~\ref{Thm:GBKnots} applied in the case $(g,b)=(t+1,0)$.

\section{High distance}
\label{Sec:HighDistance}

We now restate and prove our main theorem on high distance Heegaard
splittings for knots. 

\begin{theorem}
\label{Thm:HighDistanceKnots}
For any pair of integers $g > 1$ and  $n > 0$ there is knot $K \subset S^3$ 
and a genus $g$ splitting $S \subset E(K)$ having  distance greater than $n$.
\end{theorem}

Consider  $S^3 = V \cup_S W$ with the standard genus $g$ Heegaard
splitting.  Let $\calC_S$, $\calD_V$, and $\calD_W$ be the
corresponding curve and disk complexes.  Let $D \subset V$ be a disk
cutting $V$ into a solid torus $X$ and a handlebody $Y$ with genus
$g-1$.  Take $K_0$ to be the core of $X$.  Thus $V_0 = V - n(K_0)$ is
a compression body and $V_0 \cup W$ equals $E(K_0)$.

We must find a sequence of compression bodies $V_n \subset V$ each
homeomorphic to $V_0$ so that $\boundary_+V_n = \boundary V$ and
$$d(\calD_{V_n},\calD_W) \rightarrow \infty$$ as $n \rightarrow
\infty$.  The  knots $K_n \subset S^3$ defined by the compression
bodies $V_n$ will satisfy the
conclusion of the theorem.

When $g=2$, which is the case treated by
Johnson~\cite{Johnson06}
and Johnson-Thompson~\cite{JohnsonThompson06},
$V_0$ has a unique non-separating disk and
all other disks may be isotoped to be disjoint from it (see
e.g.~Lemma~11 in~\cite{Johnson06}).  Hence $\diam(\calD_{V_0}) = 2$
and it suffices to find a sequence $V_n$ whose non-separating disks
get arbitrarily far from $\calD_W$.  In the general case
$\diam(\calD_{V_0})=\infty$ so clearly a different construction is
needed.

To begin, set $\delta = \bdy D$.  As shown in Figure~\ref{Fig:PQ},
extend $\delta$ to a pants decomposition $\calP$ whose curves
are all meridians of $V$.  Likewise,
choose a pants decomposition $\calQ$ whose boundary curves are
meridians of $W$.  For later reference we denote by $Q_1,\ldots,Q_g$
the curves of $\calQ$ (listed in increasing size) that cross $X$ in
the front of the picture. The parallel curves on the back of the
picture are denoted by
$Q'_2,\ldots,Q'_{g-1}$;
note that $Q_1$ and $Q_g$ are their own parallels. Let $Q_{g+1}$
denote the curve going around the second hole in the picture. 

\realfig{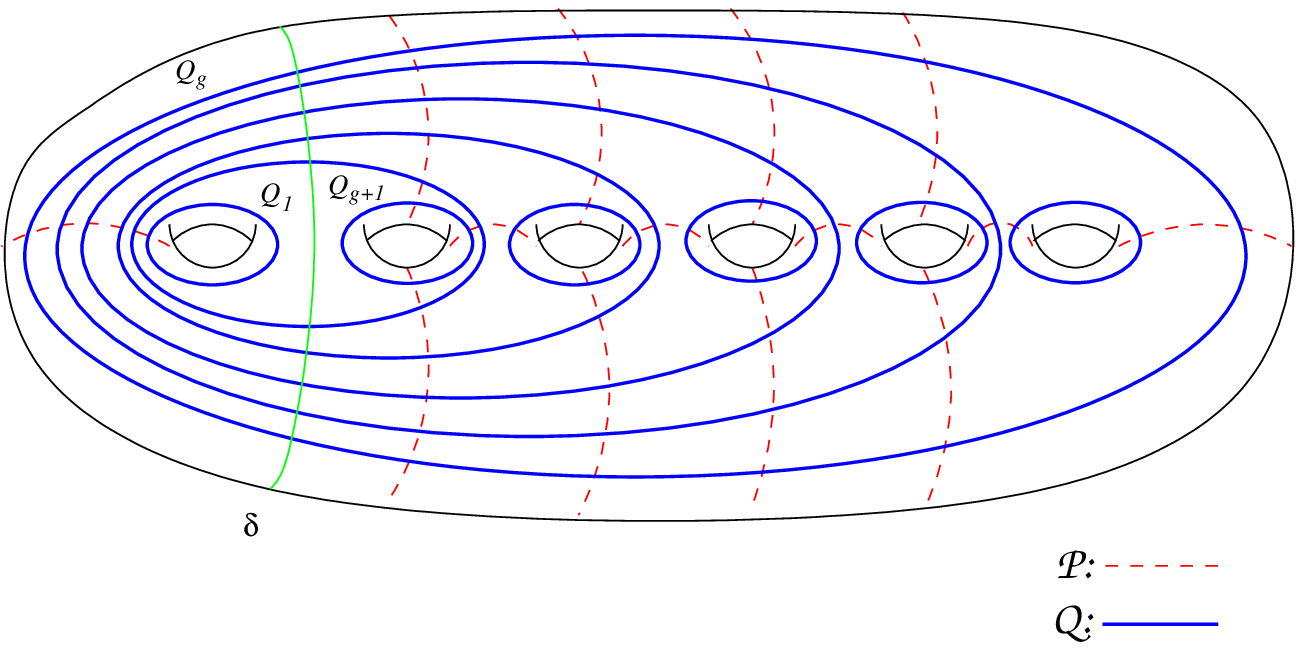}{Fig:PQ}{The pants decompositions $\calP$ of $V$ and
  $\calQ$ of $W$. Each component $Q_i$ of $\calQ$ (except the littlest and
  the biggest) has a symmetric ``partner'' $Q'_i$ on the underside of $V$.}
  
Set $S_0 = \bdy Y \cap S$.  This is the once-punctured genus $g - 1$
surface to the right of $\delta$. Let $\calP_0 = \calP \cap S_0$ be
the pair of pants decomposition of $S_0$ induced by $\calP$ (note $\calP_0$
includes $\delta$). The
important feature of our chosen decompositions is that the curves of
$\calP_0$ traverse every seam of $\calQ$ which can be isotoped rel
$\calQ$ into $S_0$.  

The construction of $V_n$ can be broken up as follows:

\subsection*{Step 1} 

Find a meridian curve $a \in \calD_V$ such that:
\begin{enumerate}
\item [$a1$:] $a$ traverses all seams of $\calQ$.
\item [$a2$:] $a$ traverses all seams of $\calP_0$.
\end{enumerate}

\subsection*{Step 2} 

Use the curve $a$ to construct a pseudo-Anosov map $\Phi \from S \to
S$ which extends over $V$ and whose stable lamination $\lambda^+$ and
unstable lamination $\lambda^-$ satisfy:
\begin{enumerate}
\item [$\lambda1$:] $\lambda^+$ traverses all seams of $\calQ$.
\item [$\lambda2$:] $\lambda^-$ traverses all seams of $\calP_0$.
\end{enumerate}

\subsection*{Step 3}  

Show that $d(\calD_W, \Phi^n(\calD_{V_0})) \to \infty$ as $n
\to \infty$.

\subsection*{Conclusion}

Since $\Phi$ extends over $V$,  we may define $V_n = \Phi^n(V_0)$
and $K_n = \Phi^n(K_0)\subset S^3$, 
so that $(V_n, W)$ is a Heegaard splitting for  the knot exterior $E(K_n)$. 
As desired in the conclusion of 
the theorem we have:
$$d(\calD_{V_n},\calD_W) \rightarrow \infty.$$

\subsection{Proof of Theorem~\ref{Thm:HighDistanceKnots}}

We now carry out the steps outlined above.

\subsection*{Step 1}  

We will first find a useful curve $\gamma \subset S_0$.  We will then
use two copies of $\gamma$ and a band sum construction to build the
desired meridian $a$.  So 
consider a train track $\tau$ adapted to $\calP_0$ and depicted in
Figure~\ref{Fig:TrainTrack}.

\realfig{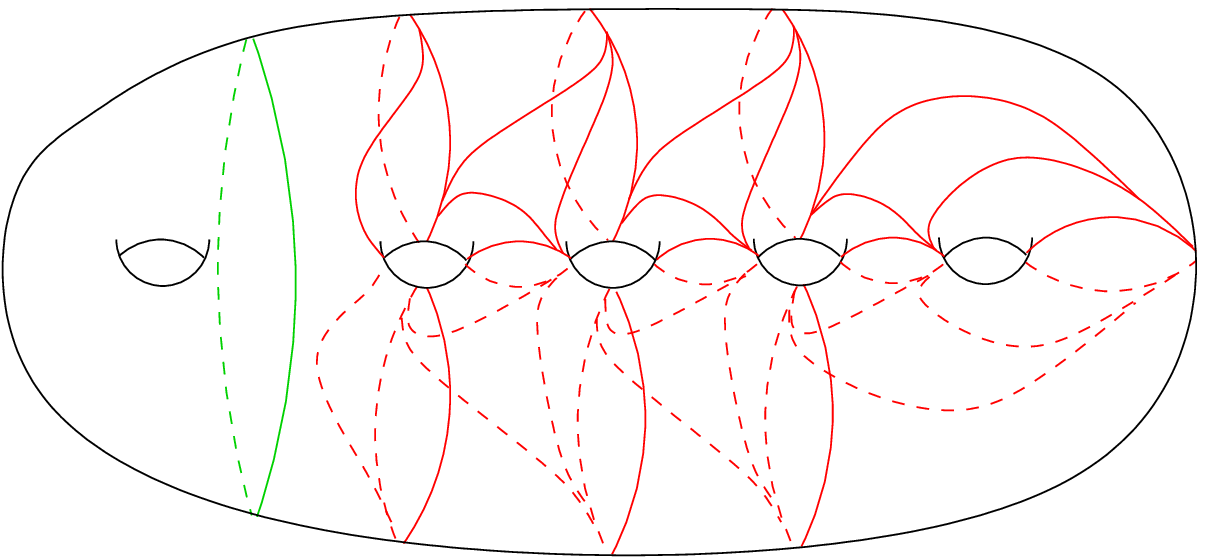}{Fig:TrainTrack}{The train track in $S_0$.}

The train track has the following properties:
\begin{enumerate}
\item 
$\tau$ contains $\calP_0 \setminus\{ \delta\}$ and additional branches. These
branches traverse all seams of $\calP_0$, with the exception of the two
seams incident to $\delta$.
\item 
The entire picture is invariant under a $180^\circ$ rotation about the
horizontal line meeting $V$ in $g+1$ arcs. 
\item The directions in which incident branches attach to each
  component of $\calP_0\setminus\{\delta\}$ on the two sides are consistent, so that
  there is a train route spiraling around any component of $\calP_0\setminus\{\delta\}$. 
\end{enumerate}
Note that these conditions determine the track up to the choice of spiraling
direction.

Observe that 
$\tau$ is a ``maximal standard train-track'' in $S_0$, in the sense of
Penner-Harer~\cite{PennerHarer92}. Hence, Theorem 3.1.2
of~\cite{PennerHarer92} tells us that positive measures on $\tau$,
modulo scaling, parameterize an open set in $\PML(S_0)$. This open set
contains a filling minimal lamination (e.g.~the stable lamination of a
pseudo-Anosov homeomorphism) since these are dense in $\PML(S_0)$.
Since simple closed curves are also dense in $\PML(S_0)$ we may
approximate this minimal lamination as well as we like by simple
closed curves carried on $\tau$. Hence there exists a simple closed
curve $\gamma$ carried on $\tau$,
traversing every branch at least twice.
Note that $\gamma$ has the following properties:
\begin{enumerate}
\item 
It traverses every seam of $\calP_0$, again excluding the two seams
incident to $\delta$.
\item 
It traverses every seam of $\calQ$ which is isotopic rel $\calQ$ into $S_0$.
\end{enumerate}
These properties follow from the corresponding properties of $\tau$,
once we observe that 
the intersections of $\tau$ with $\calQ$ are essential, i.e. that
there are no bigons in their complement, and similarly that the 
intersections of $\gamma$ with $\calP_0 \setminus \{\delta\}$
are essential.

The curve $\gamma$ fails to traverse just four seams of
$\calQ$, 
which we denote as follows: 
Let $\sigma_{1,2}$ be the seam connecting $Q_1$ to $Q_2$, in the front of $V$.
Let $\sigma'_{1,2}$ be the seam connecting $Q_1$ to $Q'_2$, in the back of $V$.
There are two seams connecting $Q_1$ and $Q_{g+1}$: one in the front, 
denoted $\sigma_{1,g+1}$, and one in the back  
denoted $\sigma'_{1,g+1}$.

\realfig{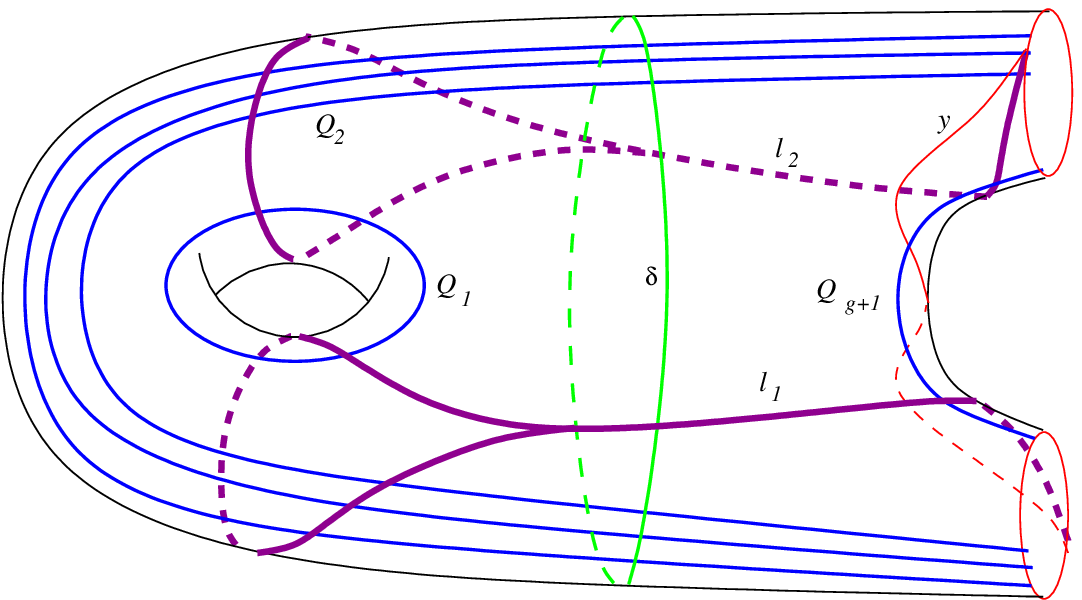}{Fig:BandSum}{The meridian $a$ is constructed via an
  enlargement $\tau'$ of $\tau$. The figure shows the 2-holed torus $Z$.
   The added loops are labeled $l_1$ and $l_2$. }

To build the meridian $a$ we add two loops $l_1$ and $l_2$ to the
train track $\tau$ 
as shown in Figure~\ref{Fig:BandSum}.  Let $Z$ denote the 2-holed
torus obtained as the union along $\delta$ of the 1-holed torus in
$\partial X$
bounded by $\delta$ 
and the adjacent pair of pants of $S_0\setminus \calP_0$.
Each loop begins inside one of
the switches of $\tau$ in $Z$, enters $X$, goes once around
the meridian disk, and then returns to the same switch.  We define an
integer measure on this new train track $\tau'$: Let $\mu$ be the
measure on $\tau$ that defines $\gamma$.  Hence $2\mu$ defines two
copies of $\gamma$.  Now subtract $2$ from the weight that $2\mu$ puts
on the branch $y$ of $\tau$ that connects the boundaries of $Z$, 
and put a weight of $1$ 
on each of the new loops.  This measure on $\tau'$ defines $a$. 
One can see that $a$ is a band sum of two copies of the meridian of
$X$ along a doubled arc of $\gamma$ -- hence it is a meridian of $V$.

Note that the intersection pattern of $\tau'$ with the curves
$Q_1,\ldots,Q_{g+1}$ and $Q'_2,\ldots,Q'_g$ has no bigons, and hence
all intersections shown are essential (the $Q'_i$ are not pictured,
but lie directly below the corresponding $Q_i$). One then checks that
$\tau'$ traverses all four seams that were previously excluded, 
namely 
$\sigma_{1,2}$,
$\sigma'_{1,2}$,
$\sigma_{1,g+1}$ and
$\sigma'_{1,g+1}$. 
Similarly there are
no bigons in the intersection pattern with $\calP_0$, and we see that
$\tau'$ traverses the two seams of $\calP_0$ incident to $\delta$. 
Since $a$ places positive measure on every branch of $\tau'$, 
it traverses all the same seams, and so satisfies conditions
$(a1)$ and $(a2)$. 

\subsection*{Step 2}  

Choose two meridians $b, c \in \calD_V$ so that $b$ and $c$ together
fill $S$. (For example:
realize $V$ as an $I$-bundle over a compact surface $F$ with exactly one
boundary component.  Let $e,f$ be
arcs in $F$ which are distance 3 or more in the arc complex of $F$ -- i.e.
there is no nontrivial arc or curve in $F$ disjoint from both $e$ and $f$ and isotopic to
neither. Then the $I$-bundles over $e$ and $f$
will be the desired disks, because
$e\union f$ cuts $F$ into disks each of which meets $\boundary F$ in a single interval. 
See also Kobayashi \cite[proof of Lemma 2.2]{Kobayashi87}.)
Let $\tau_a$, $\tau_b$ and $\tau_c$ denote the Dehn twists about
$a$, $b$ and $c$ respectively. Set $$\Phi_0 = \tau_b \circ
\tau_c^{-1}.$$ It follows from Thurston's
construction~\cite{Thurston88} that $\Phi_0$ is pseudo-Anosov, with
stable/unstable laminations $\lambda_0^\pm$. 
Since $\Phi_0$ is a composition of Dehn twists along meridian disks it
extends over $V$.  Define 
$$\Phi_N = \tau_a^N \circ \Phi_0 \circ
\tau_a^{-N}.$$ 
Since $a$ is a meridian $\Phi_N$ also extends over $V$. The stable and
unstable laminations $\lambda_N^\pm$ of $\Phi_N$ are just
$\tau_a^N(\lambda_0^\pm)$.  Since $a$ intersects $\lambda_0^\pm$ (the
latter is filling), as $N \to 
\infty$ the laminations $\lambda_N^\pm$ converge to $[a]$ in
$\mathcal{PML}(S)$.  Hence eventually both laminations satisfy
conditions $(a1)$ and $(a2)$.  Take $\Phi = \Phi_N$ for such a large
$N$ and take $\lambda^\pm = \lambda_N^\pm$.  Thus conditions $(\lambda
1)$ and $(\lambda 2)$ are satisfied.

\subsection*{Step 3} 

Since conditions $(\lambda 1)$ and $(\lambda 2)$ 
hold, we can conclude via
Lemma~\ref{Cor:FullType} that $\lambda^- \nin~\overline{\calD}_{V_0}$
and $\lambda^+ \nin~\overline{\calD}_{W}$.  
The fact that  $d(\calD_W, \Phi^n(\calD_{V_0})) \to \infty$
now follows from Lemma~\ref{Lem:Hempel}.  This completes the proof of
Theorem~\ref{Thm:HighDistanceKnots}.
\qed

\section{Ruling out simple decompositions}
\label{Sec:Decompositions}

In order to prove Theorem~\ref{Thm:GBKnots} we use a generalization by
Tomova~\cite{Tomova07} of a theorem of hers~\cite{Tomova05},
which is itself a considerable refinement of a theorem of Scharlemann
and Tomova~\cite{ScharlemannTomova05}.

\begin{theorem}
\label{Thm:Tomova}
Let $K \subset S^3$ be a knot and $S$ be a Heegaard splitting surface
of the exterior of $K$.  Suppose that $Q$ is Heegaard splitting
surface of $S^3$ so that
\begin{enumerate}
\item $K$ is in bridge position with respect to $Q$ and
\item the genus of $Q$ is less than the genus of $S$.
\end{enumerate}
Then the distance of the Heegaard splitting $S$ is at most $2 - \chi(Q
- K)$. \qed
\end{theorem}
\noindent
Condition (2) is actually stronger than necessary for Tomova's
result, but this version suffices for our needs.
We are now in position to prove Theorem~\ref{Thm:GBKnots}:

\begin{theorem}
\label{Thm:GBKnots}
For any positive integers $t$ and $b$ there is a
knot $K \subset S^3$ with tunnel number $t$ so that $K$ has no
$(t, b)$--decomposition.
\end{theorem}

\begin{proof}
Choose  a Heegaard splitting $S^3 = V \cup_S W$ of genus $t + 1$ 
and construct the knot $K \subset S^3$, as in
Theorem~\ref{Thm:HighDistanceKnots}, so that the associated splitting
$E(K) = V' \cup_S W$ has distance $d(V', W) > 2t + 2b$. 

We claim that $K$ has no $(g,c)$--decomposition for any $g\le t$, $c\le
b$. For if it did, let $Q$ be the associated splitting surface of
$S^3$. The genus $g$ of $Q$ is less than that of $S$, so we may apply
Theorem~\ref{Thm:Tomova} to obtain
$d(V',W) \leq 2 - \chi(Q - K) = 2 - (2 - 2g - 2c) \le 2t + 2b$, a
contradiction.

Setting $c=0$ we also conclude that $E(K)$ has no splitting of genus
less than $t+1$, and hence that $K$ has tunnel number $t$. 
\end{proof}

\bibliography{bibfile}
\end{document}